\def\Gn{\,\rule[-1mm]{.75mm}{4mm}\,}
\def\Z{|\hspace*{-.5mm}|\hspace*{-.5mm}|}
\def\rf#1{(\ref{#1})}
\documentstyle[12pt]{article}
\begin{document}
{\center
{\large\bf A priori bounds for Gevrey--Sobolev norms of space-periodic

\medskip
three-dimensional solutions to equations of hydrodynamic type}

\bigskip
Vladislav Zheligovsky

\medskip
International institute of earthquake prediction theory
and mathematical geophysics\\
84/32 Profsoyuznaya St., 117997 Moscow, Russian Federation\\

\medskip
UNS, CNRS, Laboratoire Cassiop\'ee, Observatoire de la C\^ote d'Azur\\
BP 4229, 06304 Nice Cedex 4, France

}

\bigskip
We present a technique for derivation of a priori bounds for Gevrey--Sobolev
norms of space-periodic three-dimensional solutions to evolutionary partial
differential equations of hydrodynamic type. It involves a transformation of
the flow velocity in the Fourier space, which introduces a feedback between
the index of the norm and the norm of the transformed solution, and results
in emergence of a mildly dissipative term. We illustrate the technique using it
to derive finite-time bounds for Gevrey--Sobolev norms of solutions to the Euler
and inviscid Burgers equations, and global in time bounds for the Voigt-type
regularisations of the Euler and Navier--Stokes equation (assuming that
the respective norm of the initial condition is bounded). The boundedness of
the norms implies analyticity of the solutions in space.

\section{Introduction}

We suggest a new approach to derivation of bounds for analytic Gevrey--Sobolev norms
(see their definition \rf{GeSo} below) of solutions to evolutionary partial
differential equations of hydrodynamic type, which exploits a feedback between
the norm of the suitably transformed solution and the index of the norm.

Whether three-dimensional solutions to the Euler equation remain smooth forever
is an open question. If the initial flow velocity belongs to the H\"older class
$C^{1+\alpha}$ (\hbox{$\alpha>0$} is arbitrary), the solution remains
in it at least for a finite time \cite{BaFr}. Moreover, if initially the flow
is analytic in space, the solution is guaranteed to be analytic in a finite
spatio-temporal region \cite{Bard1,Bard2,Bard3,Bena1,Bena2};
the loss of analyticity does not occur before the solution ceases to be
continuously differentiable \cite{Bard2}. These demonstrations involve
construction of bounds for iterated approximations of vorticity followed
along characteristics, the knowledge of the Green function in the integral
representation of the flow in the terms of vorticity, and application
of methods of the theory of functions of complex variables. It is yet unknown,
whether a singularity can develop at finite times.

Foias and Temam \cite{FoTe} proved that viscosity can enhance
the smoothness of space-periodic solutions to the Navier--Stokes equation (see
also \cite{DoTi}): If initially the flow velocity belongs to
the Sobolev space $H_1(T^3)$ (i.e., the space of space-periodic functions,
whose derivatives of order 1 are square-integrable), then at any small time
$t>0$ the flow has a bounded Gevrey--Sobolev norm. Their method was used
subsequently to demonstrate that a similar result holds for space-periodic
solutions to various modifications of the Navier--Stokes equation
\cite{KaLeTi,PaVi}. The method relies on the presence of the viscous term in
the equation, and hence is not directly applicable to the Euler equation.

We overcome the difficulty by performing a suitable transformation
of the Fourier coefficients of the flow, which introduces a dependence of
the index of a Gevrey--Sobolev norm on the norm of the transformed solution.
The transformation results in emergence of a weakly dissipative operator
in the modified equation and it becomes possible to proceed following
the method of \cite{FoTe} (time dependence of the index giving rise to a weakly
dissipative operator was employed in \cite{LeOl}). In Section 2 we use this
technique to give a simple
proof of the boundedness of a Gevrey--Sobolev norm of a space-periodic
three-dimensional solution to the Euler equation on a finite time interval.
In section 3 the same technique is applied to the Burgers equation.

The three--dimensional Navier--Stokes equation can be regularised by
introduction of suitable nonlinear \cite{La67,La68,Ladru} or linear terms, such
as hyperviscosity \cite{Lions,Sinai2}, or the time derivative of the Laplacian
of the flow (proposed by O.A.~Ladyzhenskaya at the International mathematical
congress in 1966) which can be regarded as a non-standard ``dynamic''
viscosity. Boundary value problems for the regularised equation of this type,
known as the Navier--Stokes--Voigt equation and the Oskolkov equation in the
Russian literature, were studied in \cite{Os73,Os76,Os80,Cao}. They describe
the motion of non-Newtonian visco-elastic fluids with memory (characterised
by an exponential decay of the rate of deformation under a constant stress)
arising for a class of linear integro-differential rheology equations
\cite{Os81,Os83,Os88}. While in the Voigt regularisation of the Navier--Stokes
equation dynamic viscosity is positive, its negative values may be also
physically sound \cite{Os80}. Regularisation by progressive damping
of high-wavenumber harmonics in some instances of occurrence of the flow
velocity in the equation gives rise to the so-called Camassa--Holm or
Navier--Stokes--$\alpha$\break equation \cite{FoHoTi,FoHoTi2,KaLeTi} supposedly
describing Lagrangian-averaged flows \cite{Holm}; however, an additional term
must be introduced for this interpretation to be accurate \cite{Sow2,Sow1}.

Bounds for Gevrey--Sobolev norms were derived in \cite{LaTi} for solutions to
the so-called Voigt regularisation of the Euler equation (this type of
regularisation was suggested for the Euler equation in \cite{Cao}) by the method,
initially applied in \cite{LeOl} to the ``lake equation'' (i.e., a modified
two-dimensional Euler equation). We consider a Voigt-type regularisation
of the three-dimensional Euler equation in Section 4 and show that a milder
regularisation suffices to establish the boundedness of Gevrey--Sobolev norms
of solutions. The boundedness of the norms is demonstrated for solutions
to the Voigt-type regularisation of the Navier--Stokes equation in Section 5;
the presence of diffusion enables us to further weaken the regularising term.
Constructions of Sections 4 and 5 illustrate another aspect of our approach:
dependence of the index of a Gevrey--Sobolev norm on the norm of the
transformed solution allows to decrease the order of nonlinearity in the
energy-type inequality and thus to obtain global (in time) bounds for the norm.

In this paper, we consider zero-mean (note that the means of solutions are
conserved in time) space-periodic three-dimensional solutions, assuming
for simplicity that the elementary periodicity cell is the cube
$T^3=[0,2\pi]^3$. It is not difficult to generalise our analysis to encompass
the case of arbitrary periods along the Cartesian axes. All the equations
mentioned above normally involve a force $\bf f$. Following a long-established
tradition, for the sake of simplicity we consider only the illustrative
equations for ${\bf f}=0$. Reinstating the forcing in the analysis is also
straightforward (provided the force has the right analyticity properties).

Let us also remark on the accurate procedure for derivation of a priori bounds.
The bounds are derived below for Fourier--Galerkin truncations of the solutions.
We do this by tacitly using the ODE's governing the evolution
of the Fourier coefficients of the truncated solutions. The bounds involve
norms of the initial conditions for the truncations. Assuming in their place
(in general, larger) norms of non-truncated initial conditions, we obtain bounds
for the truncated solutions, that are uniform in the number of harmonics
retained in the Fourier--Galerkin truncation. By the standard arguments, these
bounds remain valid for the norms of the non-truncated solutions.

\section{The three-dimensional Euler equation}

In this section we establish bounds for Gevrey--Sobolev norms of solutions
to the Euler equation on a finite time interval $[0,t_*)$, assuming
that the initial condition has a finite Gevrey--Sobolev norm.

The motion of a perfect (inviscid) fluid with velocity $\bf v$ under the action
of a body force $\bf f$ is governed by the Euler equation
\begin{equation}
{\partial{\bf v}\over\partial t}+({\bf v}\cdot\nabla){\bf v}={\bf f}+\nabla p.
\label{Eeq}\end{equation}
Incompressibility of fluid implies the solenoidality of the flow:
\begin{equation}
\nabla\cdot{\bf v}=0.
\label{sole}\end{equation}
The initial (at $t=0$) distribution of the flow velocity ${\bf v}^{\rm(in)}$
is prescribed.

We expand a zero-mean space-periodic flow in $T^3$ into the Fourier series:
\begin{equation}
{\bf v}=\sum_{{\bf n}\ne 0}\widehat{\bf v}_{\bf n}{\rm e}^{{\rm i}\bf n\cdot x}.
\label{Four}\end{equation}
The flow is real, as long as
\begin{equation}
\widehat{\bf v}_{\bf n}=\overline{\widehat{\bf v}_{-\bf n}},
\label{realv}\end{equation}
and solenoidality \rf{sole} implies
\begin{equation}
\widehat{\bf v}_{\bf n}\cdot{\bf n}=0.
\label{Fsole}\end{equation}

The norm $\|\cdot\|_q$ of a zero-mean space-periodic vector or scalar field
$\bf v$ in the Sobolev space $H_q(T^3)$ is defined in the terms of its Fourier
coefficients by the relation
$$\|{\bf v}\|^2_q\equiv\sum_{{\bf n}\ne0}|\widehat{\bf v}_{\bf n}|^2|{\bf n}|^{2q},$$
where $|\cdot|$ denotes the Euclidean norm of a three-dimensional vector.
For any $\sigma>0$ we define the Gevrey--Sobolev norm\footnote{Functions,
whose Gevrey--Sobolev norms are finite, are analytic (see \cite{LeOl}), and hence
the norm \rf{GeSo} might be called the analytic Gevrey--Sobolev norm. We could
extend our analysis to non-analytic Gevrey classes of index $\alpha$,
$0<\alpha<1$, whose norms involve the exponents $\exp(2\sigma|{\bf n}|^\alpha)$
instead of $\exp(2\sigma|{\bf n}|)$ in \rf{GeSo}. Under the transformation
of equations that we employ a weakly dissipative operator emerges, whose symbol
would then grow as $|{\bf n}|^\alpha$. Since the operator is the strongest
for $\alpha=1$, we consider only this case.}
\begin{equation}
\Gn{\bf w}\Gn^2_{\sigma,q}\equiv\sum_{{\bf n}\ne0}|\widehat{\bf w}_{\bf n}|^2{\rm e}^{2\sigma|{\bf n}|}|{\bf n}|^{2q}.
\label{GeSo}\end{equation}
Here the first index $\sigma$ is a lower estimate for the radius of the region
of analyticity of $\bf w$ around the real axis in the complex space. For
$\sigma=0$, \rf{GeSo} defines the norm in $H_q(T^3)$. Let $|\cdot|_q$ denote
the norm in the Lebesgue space $L_q(T^3)$. By the Sobolev embedding theorem
\cite{BeLo,Tri}, for any positive $q<3/2$ there exists a constant $C_q$ such
that for any function~$f\in H_q(T^3)$
\begin{equation}
{|f|}_{6/(3-2q)}\le C_q\|f\|_q.
\label{embe}\end{equation}

Let ${\cal P}_{\bf n}$ denote the linear operator of projection of
a three-dimensional vector on the plane normal to ${\bf n}\ne0$:
$${\cal P}_{\bf n}\widehat{\bf v}_{\bf m}\equiv\widehat{\bf v}_{\bf m}
-{\widehat{\bf v}_{\bf m}\cdot{\bf n}\over|{\bf n}|^2}{\bf n}.$$
The evolution of Fourier coefficients of the flow is governed by equations
\begin{equation}
{d\widehat{\bf v}_{\bf n}\over dt}+{\rm i}\sum_{\bf k}(\widehat{\bf v}_{\bf k}
\cdot({\bf n-k}))\,{\cal P}_{\bf n}\widehat{\bf v}_{\bf n-k}=0
\label{FEeq}\end{equation}
obtained by substitution of the series \rf{Four} (and ${\bf f}=0$) into
the Euler equation \rf{Eeq}.

Assuming $\Gn{\bf v}^{\rm(in)}\Gn_{\sigma,\,s+3/2}<\infty$ for some
positive $\sigma$ and $s\le 1/2$, we consider a transformation
\begin{equation}
\widehat{\bf v}_{\bf n}(t)=\widehat{\bf w}_{\bf n}(t)\exp\left(-\beta|{\bf n}|
\|{\bf w}({\bf x},t)\|_{s+3/2}^{-\varepsilon}\right),
\label{subs}\end{equation}
$${\bf w}({\bf x},t)=\sum_{\bf n}\widehat{\bf w}_{\bf n}(t)\,{\rm e}^{{\rm i}\bf n\cdot x},$$
where $\beta$ and $\varepsilon<2$ are positive constants. The transformation
involves solving the system of nonlinear equations \rf{subs} in
$\widehat{\bf w}_{\bf n}(t)$; the solution takes the form
$$\widehat{\bf w}_{\bf n}(t)=\widehat{\bf v}_{\bf n}(t)\exp(\psi(t)|\bf n|),$$
where the quantity $\psi(t)\ge0$ satisfies the equation
$$\psi(t)\Gn{\bf v(x},t)\Gn^\varepsilon_{\psi(t),\,s+3/2}=\beta.$$
It has a unique solution for any $t\ge0$, because the l.h.s.~is a continuous
monotonically increasing function of $\psi$ (as discussed in the introduction,
the a priori estimates are derived for Fourier--Galerkin truncations of the
flow, whereby the sum \rf{Four} is assumed to involve a finite number of terms).
If $\beta<\sigma\Gn{\bf v}^{\rm(in)}\Gn_{\sigma,\,s+3/2}^\varepsilon$, then
$\|{\bf w}({\bf x},t)\|_{s+3/2}$ is bounded at $t=0$ (uniformly over the number
of terms in the truncated sum \rf{Four}). Our goal is to show that $\bf w$
is bounded in $H_{s+3/2}(T^3)$ on some finite-length time interval $[0,t_*)$.

Substitution \rf{subs} transforms \rf{FEeq} into
$${d\widehat{\bf w}_{\bf n}\over dt}+\beta\varepsilon|{\bf n}|\,
\|{\bf w}({\bf x},t)\|_{s+3/2}^{-1-\varepsilon}\,\widehat{\bf w}_{\bf n}\,
{d\over dt}\|{\bf w}({\bf x},t)\|_{s+3/2}$$
\begin{equation}
=-{\rm i}\sum_{\bf k}(\widehat{\bf w}_{\bf k}\cdot({\bf n-k}))\,
{\cal P}_{\bf n}\widehat{\bf w}_{\bf n-k}\exp\left(\beta\|{\bf w}({\bf x},t)
\|_{s+3/2}^{-\varepsilon}\,(|{\bf n}|-|{\bf k}|-|{\bf n-k}|)\right).
\label{sEharm}\end{equation}
Scalar multiplying this equation by $|{\bf n}|^{2+2s}\,\overline{\widehat{\bf w}_{\bf n}}$
and considering the real part of the sum over $\bf n$, we find
\pagebreak
$${d\over dt}\left({1\over2}\|{\bf w}({\bf x},t)\|^2_{1+s}
+{\beta\varepsilon\over2-\varepsilon}\|{\bf w}({\bf x},t)\|^{2-\varepsilon}_{s+3/2}\right)$$
\begin{equation}
={\rm Im}\sum_{\bf k,n}(\widehat{\bf w}_{\bf k}\cdot({\bf n-k}))
(\widehat{\bf w}_{\bf n-k}\cdot\overline{\widehat{\bf w}_{\bf n}})
|{\bf n}|^{2+2s}\,\exp\left(\beta\|{\bf w}({\bf x},t)\|_{s+3/2}^{-\varepsilon}\,
(|{\bf n}|-|{\bf k}|-|{\bf n-k}|)\right).
\label{enin}\end{equation}
We call \rf{enin} the energy balance equation. In construction of estimates
for the sum arising from the nonlinear term we follow, with some variations,
the approach \cite{appl}. By the inequality for sides of
a triangle, the exponential function in the r.h.s. of \rf{enin} does not exceed
1. We will show now that the absolute value of the r.h.s.~of \rf{enin}
is bounded by $(D_s/2)\|{\bf w}\|_{s+3/2}^3$, where the constant $D_s$ is
independent of~$\bf w$ ($D_s\to\infty$ for $s\to 0$).

By virtue of the inequality
$$|{\bf n}|^{s+1/2}\le|{\bf n-k}|^{s+1/2}+|{\bf k}|^{s+1/2},$$
valid for $0<s\le 1/2$, the r.h.s.~of \rf{enin} does not exceed
$$\sum_{\bf n,k}\left(|\widehat{\bf w}_{\bf k}||{\bf n-k}|^{s+3/2}
|\widehat{\bf w}_{\bf n-k}|+|{\bf k}|^{s+1/2}|\widehat{\bf w}_{\bf k}|
|{\bf n-k}||\widehat{\bf w}_{\bf n-k}|\right)|{\bf n}|^{s+3/2}|\widehat{\bf w}_{\bf n}|$$
\begin{equation}
=(2\pi)^{-3}\int_{T^3}\left(f_0({\bf x})f_{s+3/2}({\bf x})
+f_{s+1/2}({\bf x})f_1({\bf x})\right)f_{s+3/2}(-{\bf x})\,d{\bf x},
\label{ff}\end{equation}
where scalar functions $f_q$ are defined as the Fourier series
\begin{equation}
f_q({\bf x},t)\equiv\sum_{\bf n}|\widehat{\bf w}_{\bf n}(t)||{\bf n}|^q
{\rm e}^{{\rm i}\bf n\cdot x}.
\label{fq}\end{equation}

By the Cauchy--Buniakowski--Schwarz inequality, for any
$f=\sum_{\bf n}\widehat{f}_{\bf n}{\rm e}^{{\rm i}\bf n\cdot x}$,
$$|f|\le\sum_{\bf n}|{\bf n}|^{s+3/2}|\widehat{f}_{\bf n}||{\bf n}|^{-(s+3/2)}
\le\|f\|_{s+3/2}\left(\sum_{\bf n}|{\bf n}|^{-3-2s}\right)^{1/2}.$$
The second factor in the l.h.s. of this inequality, which we denote by $c_s$,
is finite for any $s>0$ (it tends to infinity, when $s\to 0$).
Consequently, using H\"older's inequality, the Sobolev embedding theorem (see
\rf{embe}) and Parseval's identity, we find a bound for the r.h.s. of \rf{ff}:
$$(2\pi)^{-3}\left({|f_{s+3/2}|}_2\max_{T^3}|f_0|+
{|f_{s+1/2}|}_6{|f_1|}_3\right){|f_{s+3/2}|}_2$$
$$\le(2\pi)^{-3}\left(c_s\|f_0\|^2_{s+3/2}+C_1C_{1/2}\|f_0\|_{s+3/2}\|f_0\|_{3/2}\right)
\|f_0\|_{s+3/2}\le{D_s\over2}\|{\bf w}\|_{s+3/2}^3$$
for $D_s\equiv(c_s+C_1C_{1/2})/(4\pi^3)$. Thus, we obtain from \rf{enin}
$${d\over dt}\left(\|{\bf w}({\bf x},t)\|^2_{1+s}
+A\|{\bf w}({\bf x},t)\|^{2-\varepsilon}_{s+3/2}\right)
\le D_s\|{\bf w}({\bf x},t)\|_{s+3/2}^3,$$
where
\begin{equation}
A\equiv2\beta\varepsilon/(2-\varepsilon).
\label{Adef}\end{equation}
Hence, for
$\xi\equiv\|{\bf w}({\bf x},t)\|^2_{1+s}+A\|{\bf w}({\bf x},t)\|^{2-\varepsilon}_{s+3/2}$,
$${d\xi\over dt}\le D_s\,A^{-3/(2-\varepsilon)}\,\xi^{3/(2-\varepsilon)}
\quad\Rightarrow\quad-{d\over dt}\xi^{-\theta}\le D_s\theta A^{-3/(2-\varepsilon)},$$
where $\theta\equiv(1+\varepsilon)/(2-\varepsilon)$. Consequently, the bound
\begin{equation}
\xi\le\left((\xi|_{t=0})^{-\theta}-D_s\theta A^{-3/(2-\varepsilon)}\,t\right)^{-1/\theta}
\label{bnd}\end{equation}
holds for
\begin{equation}
t<t_*\equiv(D_s\theta)^{-1}A^{3/(2-\varepsilon)}(\xi|_{t=0})^{-\theta}.
\label{tstar}\end{equation}
For such $t$,
\begin{equation}
\|{\bf w}({\bf x},t)\|_{s+3/2}\le\varphi(t)\equiv
\left(\left(A/\xi|_{t=0}\right)^\theta-D_s\theta t/A\right)^{-1/(1+\varepsilon)}.
\label{wbnd}\end{equation}

Assuming in
\begin{equation}
\xi|_{t=0}=\|{\bf w}({\bf x},0)\|^2_{1+s}+A\|{\bf w}({\bf x},0)\|^{2-\varepsilon}_{s+3/2}
\label{xi0}\end{equation}
the norms of the non-truncated vector field ${\bf w}({\bf x},0)$, we obtain
bounds for $\|{\bf w}({\bf x},t)\|_{s+3/2}$ and $t_*$ that are uniform
in the number of harmonics in the Fourier--Galerkin truncation of a solution.
Inequality \rf{wbnd} and transformation \rf{subs} imply
$$\Gn{\bf v(x},t)\Gn_{\beta\varphi^{-\varepsilon},\,s+3/2}\le
\Gn{\bf v(x},t)\Gn_{\beta\|{\bf w}({\bf x},t)\|^{-\varepsilon}_{s+3/2},\,s+3/2}
=\|{\bf w}({\bf x},t)\|_{s+3/2}\le\varphi(t).$$
We have therefore proved

{\it Theorem 1.} Let initial condition ${\bf v}^{\rm(in)}$ of a solution
to the force-free Euler equation have a finite Gevrey--Sobolev norm
$\Gn{\bf v}^{\rm(in)}\Gn_{\sigma,\,s+3/2}$, where $0<s\le 1/2$. For $0\le t<t_*$
the solution satisfies the bound
\begin{equation}
\Gn{\bf v}({\bf x},t)\Gn_{\beta\varphi(t)^{-\varepsilon},\,{s+3/2}}\le\varphi(t),
\label{vbndf}\end{equation}
where $0<\beta<\sigma\Gn{\bf v}^{\rm(in)}\Gn_{\sigma,\,s+3/2}^\varepsilon$,
$0<\varepsilon<2$, $t_*$ and $\varphi(t)$ are defined by relations
\rf{tstar}--\rf{xi0}, and ${\bf w}({\bf x},0)$ is the result of application
of transformation \rf{subs} to the initial condition ${\bf v}^{\rm(in)}$.

\section{The inviscid Burgers equation}

In this section we show that our technique gives the same bound
for the Gevrey--Sobolev norms of solutions to the inviscid Burgers equation
$${\partial{\bf v}\over\partial t}+({\bf v}\cdot\nabla)\bf v=f.$$
The solenoidality of solutions is not required any more. For the sake
of simplicity, we again consider only zero-mean space-periodic solutions
for the force-free case ${\bf f}=0$, with the elementary periodicity cell
being the cube $T^3$.

A solution is expanded into the Fourier series \rf{Four};
the evolution of Fourier coefficients is now governed by equations
\begin{equation}
{d\widehat{\bf v}_{\bf n}\over dt}+{\rm i}\sum_{\bf k}(\widehat{\bf v}_{\bf k}
\cdot({\bf n-k}))\widehat{\bf v}_{\bf n-k}=0.
\label{FBeq}\end{equation}
After transformation \rf{subs} is applied, \rf{FBeq} becomes
\pagebreak
$${d\widehat{\bf w}_{\bf n}\over dt}+\beta\varepsilon|{\bf n}|\,
\|{\bf w}({\bf x},t)\|_{s+3/2}^{-1-\varepsilon}\,\widehat{\bf w}_{\bf n}\,
{d\over dt}\|{\bf w}({\bf x},t)\|_{s+3/2}$$
\begin{equation}
=-{\rm i}\sum_{\bf k}(\widehat{\bf w}_{\bf k}\cdot({\bf n-k}))\,
\widehat{\bf w}_{\bf n-k}\exp\left(\beta\|{\bf w}({\bf x},t)
\|_{s+3/2}^{-\varepsilon}\,(|{\bf n}|-|{\bf k}|-|{\bf n-k}|)\right).
\label{sBharm}\end{equation}

Scalar multiplying this equation by $|{\bf n}|^{2+2s}\,\overline{\widehat{\bf w}_{\bf n}}$
and considering the real part of the sum over $\bf n$, we again obtain
the energy balance equation \rf{enin}. Since the orthogonality \rf{Fsole}
was not used to derive from \rf{enin} the bound \rf{wbnd} yielding \rf{vbndf},
the same derivation applies for the force-free inviscid Burgers equation,
implying that the bounds \rf{wbnd} and \rf{vbndf} hold true for its solutions
as well. In other words, Theorem 1 applies literally to solutions
to the force-free inviscid Burgers equation.

\section{Voigt-type regularisation of the Euler equation}

In this section we present another illustration of our technique, demonstrating
global (in time) boundedness of Gevrey--Sobolev norms of solutions to the
Voigt-type regularisation of the Euler equation. The name ``Voigt
regularisation'' was proposed in \cite{LaTi} for $s=1$ (see \rf{Veq} below).
The boundedness was proved {\it ibid.} for $s=1$ by the method \cite{LeOl}.
Our technique seems simpler in application. We show that a milder
regularisation, than considered in \cite{LaTi}, suffices to guarantee global
regularity of solutions.

We consider three-dimensional solenoidal space-periodic zero-mean solutions
\rf{Four} to the equation
\begin{equation}
\alpha^2{\partial\over\partial t}(-\nabla^2)^s{\bf v}
+{\partial{\bf v}\over\partial t}+({\bf v}\cdot\nabla){\bf v}={\bf f}+\nabla p
\label{Veq}\end{equation}
for $s>0$. For $s=1$ this is the Voigt-regularised Euler equation, studied
in \cite{LaTi}. For ${\bf f}=0$, in the Fourier space it takes the form
\begin{equation}
(1+\alpha^2|{\bf n}|^{2s})
{d\widehat{\bf v}_{\bf n}\over dt}+{\rm i}\sum_{\bf k}(\widehat{\bf v}_{\bf k}
\cdot({\bf n-k}))\,{\cal P}_{\bf n}\widehat{\bf v}_{\bf n-k}=0.
\label{FVeq}\end{equation}

We assume $\Gn{\bf v}^{\rm(in)}\Gn_{\sigma,\,s+1/2}<\infty$ for a $\sigma>0$
and make the transformation
\begin{equation}
\widehat{\bf v}_{\bf n}(t)=\widehat{\bf w}_{\bf n}(t)\exp\left(-\beta|{\bf n}|
\,\Z{\bf w}({\bf x},t)\Z^{-\varepsilon}\right),
\label{Vsubs}\end{equation}
$${\bf w}({\bf x},t)=\sum_{\bf n}\widehat{\bf w}_{\bf n}(t)\,{\rm e}^{{\rm i}\bf n\cdot x},$$
where $\beta$ and $\varepsilon<2$ are positive constants and
\begin{equation}
\Z{\bf w}\Z^2\equiv\sum_{\bf n}(1+\alpha^2|{\bf n}|^{2s})|{\bf n}||\widehat{\bf w}_{\bf n}|^2
\label{Tnorm}\end{equation}
defines a norm equivalent to $\|\cdot\|_{s+1/2}$. Choosing
\begin{equation}
0<\beta<\sigma\left(\Gn{\bf v}^{\rm(in)}\Gn_{\sigma,1/2}^2
+\alpha^2\Gn{\bf v}^{\rm(in)}\Gn_{\sigma,\,s+1/2}^2\right)^{\varepsilon/2},
\label{betaine}\end{equation}
we ensure that
$\Z{\bf w}({\bf x},0)\Z$ is bounded uniformly over the number of Fourier
harmonics preserved in Fourier--Galerkin truncations of the initial condition
${\bf v}^{\rm(in)}$. Our goal is to derive bounds for $\Z{\bf w}({\bf x},t)\Z$.

Substitution \rf{Vsubs} transforms \rf{FVeq} into
$$(1+\alpha^2|{\bf n}|^{2s})\left({d\widehat{\bf w}_{\bf n}\over dt}
+\beta\varepsilon|{\bf n}|\,\Z{\bf w}\Z^{-1-\varepsilon}\,\widehat{\bf w}_{\bf n}\,
{d\over dt}\Z{\bf w}\Z\right)$$
$$=-{\rm i}\sum_{\bf k}(\widehat{\bf w}_{\bf k}\cdot({\bf n-k}))\,
{\cal P}_{\bf n}\widehat{\bf w}_{\bf n-k}\exp\left(\beta\Z{\bf w}\Z
^{-\varepsilon}\,(|{\bf n}|-|{\bf k}|-|{\bf n-k}|)\right).$$
The real part of the sum over $\bf n$ of these equations scalar multiplied
by $\overline{\widehat{\bf w}_{\bf n}}$ reduces to
$${d\over dt}\left(\|{\bf w}\|^2_0+\alpha^2\|{\bf w}\|^2_s+A\Z{\bf w}\Z^{2-\varepsilon}\right)$$
$$=2{\rm Im}\sum_{\bf k,n}(\widehat{\bf w}_{\bf k}\cdot({\bf n-k}))
(\widehat{\bf w}_{\bf n-k}\cdot\overline{\widehat{\bf w}_{\bf n}})
\exp\left(\beta\Z{\bf w}\Z^{-\varepsilon}\,(|{\bf n}|-|{\bf k}|-|{\bf n-k}|)\right)$$
$$={\rm Im}\sum_{\bf k,n}(\widehat{\bf w}_{\bf k}\cdot({\bf n-k}))
(\widehat{\bf w}_{\bf n-k}\cdot\overline{\widehat{\bf w}_{\bf n}})
\left(\exp\left(\beta\Z{\bf w}\Z^{-\varepsilon}\,(|{\bf n}|-|{\bf k}|-|{\bf n-k}|)\right)\right.$$
\begin{equation}
-\left.\exp\left(\beta\Z{\bf w}\Z^{-\varepsilon}\,(|{\bf n-k}|-|{\bf k}|-|{\bf n}|)\right)\right)
\label{Venin}\end{equation}
(above, in one of the sums in $\bf k$ and $\bf n$ we have changed the index
of summation $\bf n\to k-n$ and relied on the fact that $\bf w$ is real-valued
and solenoidal, which stems from \rf{Fsole} and \rf{realv}).

Note that $|e^{\alpha'}-e^{\alpha''}|\le|\alpha'-\alpha''|$ for any negative
$\alpha'$ and $\alpha''$. Consequently, the absolute value of the r.h.s.~of
\rf{Venin} does not exceed
$$\sum_{\bf n,k}|\widehat{\bf w}_{\bf k}||{\bf n-k}||\widehat{\bf w}_{\bf n-k}|
|\widehat{\bf w}_{\bf n}|\,2\beta\,\Z{\bf w}\Z^{-\varepsilon}\,||{\bf n}|-|{\bf n-k}||
\le{\beta\Z{\bf w}\Z^{-\varepsilon}\over4\pi^3}\int_{T^3}f_0(-{\bf x})f_1^2({\bf x})\,d{\bf x}$$
(here and in \rf{Venin} we use the notation $A$ and $f_q$ \rf{fq} introduced in Section 2)
\begin{equation}
\le{\beta\Z{\bf w}\Z^{-\varepsilon}\over4\pi^3}|f_0|_{9/(2-3\zeta)}|f_1|_{18/(7+3\zeta)}^2
\label{Vff}\end{equation}
(by H\"older's inequality for a sufficiently small $\zeta>0$).
By the Sobolev embedding theorem
$$|f_0|_{9/(2-3\zeta)}\le C_{5/6+\zeta}\|f_0\|_{5/6+\zeta},\qquad
|f_1|_{18/(7+3\zeta)}\le C_{1/3-\zeta/2}\|f_1\|_{1/3-\zeta/2}$$
(see \rf{embe}).
By virtue of \rf{fq} and by H\"older's inequality, \rf{Vff} is bounded by
$${\beta C_{5/6+\zeta}C_{1/3-\zeta/2}^2\Z{\bf w}\Z^{-\varepsilon}\over4\pi^3}
\|{\bf w}\|_{5/6+\zeta}\|{\bf w}\|^2_{4/3-\zeta/2}
\le{\beta C_{5/6+\zeta}C_{1/3-\zeta/2}^2\Z{\bf w}\Z^{-\varepsilon}\over4\pi^3}
\|{\bf w}\|^{1+6\zeta}_{5/6+\zeta}\|{\bf w}\|^{2-6\zeta}_{4/3+\zeta}.$$
Hence, for $s\ge5/6+\zeta$ and $\varepsilon=2-6\zeta$, provided $0<\zeta<1/3$,
we obtain from the energy balance equation \rf{Venin}:
$${d\over dt}\left(\|{\bf w}\|^2_0+\alpha^2\|{\bf w}\|^2_s+A\Z{\bf w}\Z^{6\zeta}\right)
\le{\beta C_{5/6+\zeta}C_{1/3-\zeta/2}^2\over4\pi^3\alpha^\varepsilon}\|{\bf w}\|^{1+6\zeta}_s.$$
This implies
\begin{equation}
{d\xi\over dt}\le D\,\xi^{1/2+3\zeta},
\label{Vxiin}\end{equation}
where it is denoted
$$D\equiv\beta C_{5/6+\zeta}C_{1/3-\zeta/2}^2(\pi\alpha)^{-3}/4,$$
\begin{equation}
\xi\equiv\|{\bf w}\|^2_0+\alpha^2\|{\bf w}\|^2_s+A\Z{\bf w}\Z^{6\zeta}.
\label{xidef}\end{equation}
For $0<\zeta<1/6$, we obtain from \rf{Vxiin} a global polynomial bound
\begin{equation}
\xi\le\left((\xi|_{t=0})^{1/2-3\zeta}+(1/2-3\zeta)\,Dt\right)^{(1/2-3\zeta)^{-1}}
\equiv\varphi(t;\zeta),
\label{Vpb}\end{equation}
and for $\zeta=1/6$ an exponential one
\begin{equation}
\xi\le\xi|_{t=0}\,{\rm e}^{Dt}\equiv\varphi(t;1/6).
\label{Veb}\end{equation}
For $\zeta>1/6$, \rf{Vxiin} yields bounds for $\xi$, which are finite-time
and thus not of interest. In view of \rf{Vsubs}, \rf{Vpb} and
\rf{Veb}, the solution to equation \rf{Veq} satisfies the inequalities
$$\Z{\bf w}\Z\le(\varphi/A)^{1/(6\zeta)}\qquad\Rightarrow$$
$$\Gn{\bf v}({\bf x},t)\Gn_{\beta(\varphi(t;\zeta)/A)^{1-1/(3\zeta)},\,1/2}\le
\Gn{\bf v}({\bf x},t)\Gn_{\beta\Z{\bf w}({\bf x},t)\Z^{-\varepsilon},\,1/2}=
\|{\bf w}\|_{1/2}\le\Z{\bf w}\Z\le(\varphi(t;\zeta)/A)^{1/(6\zeta)}.$$
Note that the first index of the norm in the l.h.s.~of this inequality
is strictly positive for any $t>0$.

We have proved the following statement:

{\it Theorem 2.} Suppose $0<\zeta\le1/6$ and $s\ge5/6+\zeta$. Let the norm
$\Gn{\bf v}^{\rm(in)}\Gn_{\sigma,\,s+1/2}$ of initial condition ${\bf v}^{\rm(in)}$
of a solution to the Voigt-type regularisation of the force-free
three-dimensional Euler equation \rf{Veq} be finite for some $\sigma>0$.
Then the bound
$$\Gn{\bf v}({\bf x},t)\Gn_{\beta(\varphi(t;\zeta)/A)^{1-1/(3\zeta)},\,1/2}\le(\varphi(t;\zeta)/A)^{1/(6\zeta)}$$
holds true for the solution at any time $t>0$. Here $A$ is defined by \rf{Adef},
$\varepsilon=2-6\zeta$ and $\beta$ satisfies inequality \rf{betaine};
$\varphi$ is defined by \rf{Vpb} and \rf{Veb}, where $\xi|_{t=0}$ is obtained
by application of formula \rf{xidef} to ${\bf w}({\bf x},0)$, which is
the result of application of transformation \rf{subs} to the non-truncated
initial condition ${\bf v}^{\rm(in)}$.

To the best of our knowledge, the regularised equation \rf{Veq} with
the fractional power of the Laplacian has not yet been considered in Sobolev
spaces in literature. It is easy to show incrementally that for $s\ge5/6$
the solution of \rf{Veq} belongs to Sobolev spaces of arbitrary high indices
(provided the initial data and the forcing $\bf f$ are sufficiently regular).
Such a proof can be started by multiplication of \rf{Veq} by $\bf v$ and
deriving a bound for $\|{\bf v}\|_{5/6}$. Scalar multiplying now \rf{Veq}
by $(-\nabla^2)^{5/3-s}\bf v$, using the inequality
$|{\bf v}|_{9/2}\le C_{5/6}\|{\bf v}\|_{5/6}$ to bound the integral arising
from the nonlinear term, it easy to show that $\|{\bf v}\|_{5/3}$ is bounded
at any $t>0$. Higher-index Sobolev norms can be subsequently bounded similarly.

\section{Voigt-type regularisation\\ of the Navier--Stokes equation}

In this section we prove global (in time) boundedness of Gevrey--Sobolev norms
of three-dimensional solenoidal space-periodic zero-mean solutions \rf{Four}
to the Voigt-type regularisation of the Navier--Stokes equation
\begin{equation}
\alpha^2{\partial\over\partial t}(-\nabla^2)^s{\bf v}+{\partial{\bf v}\over\partial t}
+({\bf v}\cdot\nabla){\bf v}-\nu\nabla^2{\bf v}={\bf f}+\nabla p,
\label{VNSeq}\end{equation}
where $s>0$. The Voigt regularisation ($s=1$) was previously investigated in
\cite{Os73,Zh}; boundedness of Gevrey--Sobolev norms of solutions for $s=1$ was
demonstrated in \cite{KaLeTi} by a method, different from the one used here.
We show here that a milder regularisation, than the one considered
in the previous section for the Euler equation, suffices to guarantee global
regularity of solutions. In the first two subsections we derive the bounds
separately for $s>1/2$, where the smoothing term is excessively strong, and for
the limit value $s=1/2$ --- constructions in the two cases are
somewhat different. In the last subsection we show that the method \cite{FoTe}
can be used for $s\le1/2$ to establish instantaneous development of analyticity
of solutions for mildly regular initial conditions.

\subsection{Excessive damping: $s>1/2$}

Assuming $\Gn{\bf v}^{\rm(in)}\Gn_{\sigma,\,s+1/2}<\infty$ for a positive
$\sigma$, we make transformation \rf{Vsubs} for $1<\varepsilon<2$, where
the norm $\Z\cdot\Z$, equivalent to $\|\cdot\|_{s+1/2}$, is defined
by \rf{Tnorm}; $\beta$ is supposed to satisfy inequality \rf{betaine},
implying $\Z{\bf w}({\bf x},0)\Z<\infty$. We need to derive a bound
for $\Z{\bf w}({\bf x},t)\Z$.

Substituting \rf{Four} and \rf{Vsubs}, we transform \rf{VNSeq} (for
${\bf f}=0$) into the system of equations
$$(1+\alpha^2|{\bf n}|^{2s})\left({d\widehat{\bf w}_{\bf n}\over dt}
+\beta\varepsilon|{\bf n}|\,\Z{\bf w}\Z^{-1-\varepsilon}\,\widehat{\bf w}_{\bf n}\,
{d\over dt}\Z{\bf w}\Z\right)+\nu|{\bf n}|^2\widehat{\bf w}_{\bf n}$$
$$=-{\rm i}\sum_{\bf k}(\widehat{\bf w}_{\bf k}\cdot({\bf n-k}))\,
{\cal P}_{\bf n}\widehat{\bf w}_{\bf n-k}\exp\left(\beta\Z{\bf w}\Z
^{-\varepsilon}\,(|{\bf n}|-|{\bf k}|-|{\bf n-k}|)\right).$$
The real part of the sum of these equations scalar multiplied by
$\overline{\widehat{\bf w}_{\bf n}}$ reduces to
$${d\over dt}\left(\|{\bf w}\|^2_0+\alpha^2\|{\bf w}\|^2_s
+A\Z{\bf w}\Z^{2-\varepsilon}\right)+2\nu\|{\bf w}\|_1^2$$
$$={\rm Im}\sum_{\bf k,n}(\widehat{\bf w}_{\bf k}\cdot({\bf n-k}))
(\widehat{\bf w}_{\bf n-k}\cdot\overline{\widehat{\bf w}_{\bf n}})
\left(\exp\left(\beta\Z{\bf w}\Z^{-\varepsilon}\,(|{\bf n}|-|{\bf k}|-|{\bf n-k}|)\right)\right.$$
\begin{equation}
-\left.\exp\left(\beta\Z{\bf w}\Z^{-\varepsilon}\,(|{\bf n-k}|-|{\bf k}|-|{\bf n}|)\right)\right)
\label{VNSenin}\end{equation}
(this identity differs from \rf{Venin} only by the presence of
the diffusion-related term $2\nu\|{\bf w}\|_1^2$).

Note that $|e^{\alpha'}-e^{\alpha''}|\le|\alpha'-\alpha''|^{1/\varepsilon}$ for
all $\alpha'\le0$ and $\alpha''\le0$ (since $\varepsilon>1$). Consequently,
the r.h.s.~of \rf{VNSenin} does not exceed
\pagebreak
$$\sum_{\bf n,k}|\widehat{\bf w}_{\bf k}||{\bf n-k}||\widehat{\bf w}_{\bf n-k}|
|\widehat{\bf w}_{\bf n}|\,(2\beta\Z{\bf w}\Z^{-\varepsilon}\,||{\bf n}|-|{\bf n-k}||)^{1/\varepsilon}$$
$$\le{(2\beta)^{1/\varepsilon}\over(2\pi)^3\Z{\bf w}\Z}\int_{T^3}f_0(-{\bf x})
f_{1/\varepsilon}({\bf x})f_1({\bf x})\,d{\bf x}
\le{(2\beta)^{1/\varepsilon}\over(2\pi)^3\Z{\bf w}\Z}
{|f_0|}_6{|f_{1/\varepsilon}|}_{3/(s+1/2)}{|f_1|}_{3/(2-s)}$$
(by H\"older's inequality and definition \rf{fq} of the functions $f_q$)
\begin{equation}
\le{(2\beta)^{1/\varepsilon}\over(2\pi)^3\Z{\bf w}\Z}
C_1{\|f_0\|}_1C_{1-s}{\|f_{1/\varepsilon}\|}_{1-s}C_{s-1/2}{\|f_1\|}_{s-1/2}
\label{VNSff}\end{equation}
(by the Sobolev embedding theorem, see \rf{embe}). We have assumed here
$1/2<s<1$, and we further demand $2>\varepsilon>1/s$. Then, by H\"older's
inequality and definition \rf{fq} of the functions $f_q$, \rf{VNSff} is
bounded by
$${C_1C_{1-s}C_{s-1/2}(2\beta)^{1/\varepsilon}\over(2\pi)^3\Z{\bf w}\Z}
\|{\bf w}\|_1\|{\bf w}\|_{1-s+1/\varepsilon}\|{\bf w}\|_{s+1/2}
\le D_{s,\varepsilon}\|{\bf w}\|_1^{2-\kappa}\|{\bf w}\|_s^\kappa,$$
where we have denoted
$$D_{s,\varepsilon}\equiv C_1C_{1-s}C_{s-1/2}(2\beta)^{1/\varepsilon}(2\pi)^{-3}/\alpha,
\qquad\kappa\equiv\min(1,(s-1/\varepsilon)/(1-s)).$$

We therefore obtain from \rf{VNSenin}, by Young's inequality,
$${d\over dt}\left(\|{\bf w}\|^2_0+\alpha^2\|{\bf w}\|^2_s+A\Z{\bf w}\Z^{2-\varepsilon}\right)
\le D'_{s,\varepsilon}\alpha^2\|{\bf w}\|_s^2,$$
where
$$D'_{s,\varepsilon}\equiv{D^{2/\kappa}_{s,\varepsilon}\kappa\over2\alpha^2}
\left({2-\kappa\over4\nu}\right)^{(2-\kappa)/\kappa}.$$
Therefore,
\begin{equation}
\|{\bf w}\|^2_0+\alpha^2\|{\bf w}\|^2_s+A\Z{\bf w}\Z^{2-\varepsilon}\le\left.
\left(\|{\bf w}\|^2_0+\alpha^2\|{\bf w}\|^2_s+A\Z{\bf w}\Z^{2-\varepsilon}\right)\right|_{t=0}
{\rm e}^{D'_{s,\varepsilon}t}\equiv\varphi(t).
\label{phidef}\end{equation}
Since $\Z{\bf w}\Z\le(\varphi/A)^{1/(2-\varepsilon)}$, transformation
\rf{Vsubs} implies
$$\Gn{\bf v}({\bf x},t)\Gn_{\beta(\varphi(t)/A)^{-\varepsilon/(2-\varepsilon)},\,1/2}\le
\Gn{\bf v}({\bf x},t)\Gn_{\beta\Z{\bf w}({\bf x},t)\Z^{-\varepsilon},\,1/2}=
\|{\bf w}\|_{1/2}\le\Z{\bf w}\Z\le(\varphi(t)/A)^{1/(2-\varepsilon)}.$$
The first index of the Gevrey--Sobolev norm in the l.h.s.~of this inequality is
strictly positive for any $t>0$.

Thus, we have demonstrated

{\it Theorem 3.} Suppose $2>\varepsilon>1/s>1$. Let the norm
$\Gn{\bf v}^{\rm(in)}\Gn_{\sigma,\,s+1/2}$ of initial condition ${\bf v}^{\rm(in)}$
of a solution to the Voigt-type regularisation of the force-free
three-dimensional Navier--Stokes equation \rf{VNSeq} be finite for a $\sigma>0$.
Then the bound
$$\Gn{\bf v}({\bf x},t)\Gn_{\beta(\varphi(t)/A)^{-\varepsilon/(2-\varepsilon)},\,1/2}
\le(\varphi(t)/A)^{1/(2-\varepsilon)}$$
holds true for the solution at any time $t>0$. Here $A$ is defined by
\rf{Adef}, $\beta$ satisfies inequality \rf{betaine} and $\varphi$ is defined
by \rf{phidef}, where ${\bf w}({\bf x},0)$ is obtained from the non-truncated
initial condition ${\bf v}^{\rm(in)}$ by transformation \rf{subs}.

\subsection{The critical damping: $s=1/2$}

In this subsection we focus on the limit value $s=1/2$ and show that in
this case bounds for the Gevrey--Sobolev norms can be obtained in a similar way.
Instead of \rf{Vsubs}, we now make a substitution
\begin{equation}
\widehat{\bf v}_{\bf n}(t)=\widehat{\bf w}_{\bf n}(t)\exp\left(-\beta|{\bf n}|
(1+\Z{\bf w}({\bf x},t)\Z)^{-2}\right).
\label{hVsubs}\end{equation}
$\Z{\bf w}({\bf x},0)\Z$ is bounded uniformly
over truncations of ${\bf v}^{\rm(in)}$ provided
\begin{equation}
0<\beta<\sigma\left(1+\sqrt{\Gn{\bf v}^{\rm(in)}\Gn_{\sigma,1/2}^2
+\alpha^2\Gn{\bf v}^{\rm(in)}\Gn_{\sigma,1}^2}\ \right)^2
\label{newbeta}\end{equation}
(under this condition a transformation of the non-truncated initial condition
${\bf v}^{\rm(in)}$ is also well-defined).

As in the previous subsection, we derive from the Voigt-type regularised
Navier--Stokes equation \rf{VNSeq} (for ${\bf f}=0$) a system of equations
governing the evolution of Fourier coefficients $\widehat{\bf w}_{\bf n}(t)$
of $\bf w$, and consider the real part of the sum of these equations, scalar
multiplied by $\overline{\widehat{\bf w}_{\bf n}}$. The energy
balance equation, analogous to \rf{VNSenin}, takes now the form
$${d\over dt}(\|{\bf w}\|^2_0+\alpha^2\|{\bf w}\|^2_{1/2})+4\beta\,
{\Z{\bf w}\Z^2\over(1+\Z{\bf w}\Z)^3}{d\over dt}\Z{\bf w}\Z+2\nu\|{\bf w}\|_1^2$$
$$={\rm Im}\sum_{\bf k,n}(\widehat{\bf w}_{\bf k}\cdot({\bf n-k}))
(\widehat{\bf w}_{\bf n-k}\cdot\overline{\widehat{\bf w}_{\bf n}})
\left(\exp\left(\beta(|{\bf n}|-|{\bf k}|-|{\bf n-k}|)/(1+\Z{\bf w}\Z)^2\right)\right.$$
\begin{equation}
-\left.\exp\left(\beta(|{\bf n-k}|-|{\bf k}|-|{\bf n}|)/(1+\Z{\bf w}\Z)^2\right)\right).
\label{VNSbndh}\end{equation}
Using the inequality $|e^{\alpha'}-e^{\alpha''}|\le|\alpha'-\alpha''|^{1/2}$
which holds for any $\alpha'\le0$ and $\alpha''\le0$,
we obtain from \rf{VNSbndh} the energy inequality
$${d\over dt}\left(\|{\bf w}\|^2_0+\alpha^2\|{\bf w}\|^2_{1/2}+4\beta\,\left(
\ln(1+\Z{\bf w}\Z)+{4\Z{\bf w}\Z+3\over2(1+\Z{\bf w}\Z)^2}\right)\right)+2\nu\|{\bf w}\|_1^2$$
$$\le\sqrt{2\beta}\sum_{\bf n,k}|\widehat{\bf w}_{\bf k}||{\bf n-k}||\widehat{\bf w}_{\bf n-k}|
|\widehat{\bf w}_{\bf n}||{\bf k}|^{1/2}(1+\Z{\bf w}\Z)^{-1}
\le{\sqrt{2\beta}\over(2\pi)^3}C_{1/2}C_1{\|{\bf w}\|_1^3\over1+\Z{\bf w}\Z}$$
\begin{equation}
\le{\sqrt{2\beta}C_{1/2}C_1\over(2\pi)^3\alpha^2}\|{\bf w}\|_1^2.
\label{bndh}\end{equation}
Consequently, if
$$\beta\le2\nu^2(2\pi)^6C^{-2}_{1/2}C^{-2}_1\alpha^4$$
(in addition to condition \rf{newbeta}), then for
\begin{equation}
\xi(t)\equiv\|{\bf w}\|^2_0+\alpha^2\|{\bf w}\|^2_{1/2}+4\beta\ln(1+\Z{\bf w}\Z)
\label{xiln}\end{equation}
we find, integrating \rf{bndh} in time:
$$\xi(t)\le\xi(0)+4\beta\left(-{4\Z{\bf w}({\bf x},t)\Z+3\over2(1+\Z{\bf w}({\bf x},t)\Z)^2}
+{4\Z{\bf w}({\bf x},0)\Z+3\over2(1+\Z{\bf w}({\bf x},0)\Z)^2}\right)
\le\xi(0)+6\beta.$$
Each of the three terms, constituting $\xi(t)$, is positive; hence we deduce
from this inequality
$$1+\Z{\bf w}({\bf x},t)\Z\le\exp(\xi(0)/(4\beta)+3/2),\qquad
\|{\bf w}({\bf x},t)\|_0\le\sqrt{\xi(0)+6\beta}.$$
In view of \rf{hVsubs}, the two inequalities imply
$$\Gn{\bf v}({\bf x},t)\Gn_{\beta\exp(-\xi(0)/(2\beta)-3),\,0}\le
\Gn{\bf v}({\bf x},t)\Gn_{\beta(1+\Z{\bf w}({\bf x},t)\Z)^{-2},\,0}\le
\|{\bf w}({\bf x},t)\|_0\le\sqrt{\xi(0)+6\beta}.$$
Thus, in the case of critical damping the following theorem holds:

{\it Theorem 4.} Suppose $\sigma>0$ and the norm
$\Gn{\bf v}^{\rm(in)}\Gn_{\sigma,1}$ of initial condition ${\bf v}^{\rm(in)}$
of a solution to the Voigt-type regularisation of the force-free
three-dimensional Navier--Stokes equation \rf{VNSeq} for $s=1/2$ is finite.
Then the bound
\begin{equation}
\Gn{\bf v}({\bf x},t)\Gn_{\beta\exp(-\xi(0)/(2\beta)-3),\,0}\le\sqrt{\xi(0)+6\beta}
\label{hNSVGS}\end{equation}
holds true for the solution at any time $t>0$. Here $\beta$ satisfies
inequality \rf{newbeta} and $\xi(0)$ is determined by application of formula
\rf{xiln} to ${\bf w}({\bf x},0)$, obtained from the non-truncated initial
condition ${\bf v}^{\rm(in)}$ by transformation \rf{hVsubs}.

What happens if the initial data ${\bf v}^{\rm(in)}$ is non-analytic and thus
Theorem 4 is inapplicable? Like in the case of the Voigt-type regularisation
of the three-dimensional Euler equation, it is possible to develop the theory
of solutions to equation \rf{VNSeq} in the Sobolev spaces. In particular, one
can show incrementally that for $s\ge1/2$ the solution of \rf{VNSeq} belongs
to Sobolev spaces of arbitrary high indices (when the initial data and
the forcing $\bf f$ are sufficiently regular). Since this question is not
in the scope of our paper, we only present a brief sketch of derivation
of the bounds. Multiplication of \rf{VNSeq} by $\bf v$ demonstrates boundedness
of $\int_0^t\|{\bf v(x},\tau)\|_1^2d\tau$ for any $t>0$. Scalar multiplication
of the equation by $(-\nabla^2)\bf v$ and the use of the inequality
$$|({\bf v}\cdot\nabla){\bf v}\cdot(-\nabla^2){\bf v}|\le
|{\bf v}|_6|\nabla{\bf v}|_3|\nabla^2{\bf v}|_2\le
(C\|{\bf v}\|^2_1\|{\bf v}\|^2_{1+s}+\nu\|{\bf v}\|_2^2)/2$$
yields
$${d\over dt}(\alpha^2\|{\bf v}\|^2_{1+s}+\|{\bf v}\|_1^2)\le
\|{\bf f}\|^2_{1-s}+(1+C\|{\bf v}\|^2_1)\|{\bf v}\|^2_{1+s},$$
whereby
$$\alpha^2\|{\bf v}\|^2_{1+s}+\|{\bf v}\|_1^2\le
(\alpha^2\|{\bf v(x},0)\|^2_{1+s}+\|{\bf v(x},0)\|_1^2)\exp\left(\alpha^{-2}
\int_0^t(1+C\|{\bf v(x},\tau)\|^2_1)d\tau\right)$$
$$+\int_0^t\|{\bf f(x},\tau)\|^2_{1-s}\exp\left(\alpha^{-2}
\int_\tau^t(1+C\|{\bf v(x},\tau')\|^2_1)d\tau'\right)d\tau.$$
Bounds for higher-index Sobolev norms can be subsequently derived by a similar
procedure.

\subsection{Instantaneous development of analyticity for $s\le1/2$}

The remark, made in \cite{KaLeTi}, that the Navier--Stokes--Voigt equation
($s=1$) exhibits damped hyperbolicity, remains valid for the milder
regularisation for $s>1/2$ considered in subsection 5.1. In particular,
solutions for $s>1/2$ apparently cannot instantaneously acquire analyticity ---
at least, the method \cite{FoTe}, revealing that solutions to
the Navier--Stokes equation are capable of this, is not directly applicable
to equation \rf{VNSeq}. To see this, let us transform a solution to \rf{VNSeq},
following \cite{FoTe}, using the relations
\begin{equation}
{\bf v}({\bf x},t)=\sum_{{\bf n}\ne 0}\widetilde{\bf w}_{\bf n}(t)
{\rm e}^{-\beta t|{\bf n}|+{\rm i}\bf n\cdot x},\qquad
{\bf w}({\bf x},t)=\sum_{\bf n}\widetilde{\bf w}_{\bf n}(t)\,{\rm e}^{{\rm i}\bf n\cdot x},
\label{TFtr}\end{equation}
where $\beta>0$ is a constant. The evolution of the transformed Fourier
coefficients satisfies the equation
$$(1+\alpha^2|{\bf n}|^{2s})\left({d\widetilde{\bf w}_{\bf n}\over dt}
-\beta|{\bf n}|\widetilde{\bf w}_{\bf n}\right)+\nu|{\bf n}|^2\widetilde{\bf w}_{\bf n}$$
$$=-{\rm i}\sum_{\bf k}(\widetilde{\bf w}_{\bf k}\cdot({\bf n-k}))\,
{\cal P}_{\bf n}\widetilde{\bf w}_{\bf n-k}\exp\left(\beta t(|{\bf n}|-|{\bf k}|-|{\bf n-k}|)\right).$$
Scalar multiplying it by $|{\bf n}|^{2\gamma}\,\overline{\widetilde{\bf w}_{\bf n}}$
for $\gamma\ge0$, summing over $\bf n$ and taking the real part we find
$${1\over2}\,{d\over dt}(\|{\bf w}\|^2_\gamma+\alpha^2\|{\bf w}\|^2_{\gamma+s})
-\beta(\|{\bf w}\|^2_{1/2+\gamma}+\alpha^2\|{\bf w}\|^2_{1/2+\gamma+s})+\nu\|{\bf w}\|_{1+\gamma}^2$$
\begin{equation}
={\rm Im}\sum_{\bf k,n}(\widetilde{\bf w}_{\bf k}\cdot({\bf n-k}))
(\widetilde{\bf w}_{\bf n-k}\cdot\overline{\widetilde{\bf w}_{\bf n}})
|{\bf n}|^{2\gamma}\exp\left(\beta t(|{\bf n}|-|{\bf k}|-|{\bf n-k}|)\right).
\label{last}\end{equation}
Thus, the strength of the regularising term for $s>1/2$
gives rise to a problem: viscous dissipation is too weak
to control the term $\beta\alpha^2\|{\bf w}\|^2_{1/2+\gamma+s}$, appearing
in the l.h.s.~of the energy balance equation \rf{last}. Still, for $s>1/2$
finite-time Gevrey class $G_{1/(2-2s)}$ regularity emerges instantaneously,
and this can be established by the method \cite{FoTe}.

To show that the method \cite{FoTe} works for $0<s\le1/2$, we choose three
quantities $\eta_i>0$ such that $\eta_1+\eta_2+\eta_3=\nu$ and
$1/2<\gamma\le1$. If $s=1/2$, we also demand
\begin{equation}
\beta\le\eta_2/\alpha^2
\label{xbeta}\end{equation}
in transformation \rf{TFtr}. By H\"older's and Young's inequalities,
\begin{equation}
\beta\|{\bf w}\|^2_{1/2+\gamma}\le\eta_1\|{\bf w}\|_{1+\gamma}^2+
(4\eta_1)^{-1}\beta^2\|{\bf w}\|^2_\gamma,
\label{t1}\end{equation}
\begin{equation}
\beta\alpha^2\|{\bf w}\|^2_{1/2+\gamma+s}\le\eta_2\|{\bf w}\|_{1+\gamma}^2
+Q_1\alpha^2\|{\bf w}\|^2_{\gamma+s},
\label{t2}\end{equation}
where
$$Q_1\equiv\left\{
\begin{array}{ll}
0,&\mbox{\ if\ }\beta\le\eta_2/\alpha^2,\\
\displaystyle{\beta(1-2s)\over2(1-s)}\left({\beta\alpha^2\over2\eta_2(1-s)}\right)^{1/(1-2s)},
&\mbox{\ otherwise.}
\end{array}\right.$$
Using H\"older's inequality, the Sobolev embedding theorem \rf{embe}
and Young's inequality, we find that the r.h.s. of \rf{last} is bounded by the sum
\pagebreak
$$\sum_{\bf k,n}|\widetilde{\bf w}_{\bf k}||{\bf n-k}|
|\widetilde{\bf w}_{\bf n-k}||\widetilde{\bf w}_{\bf n}||{\bf n}|^{2\gamma}$$
$$=(2\pi)^{-3}\int_{T^3}
\left(\sum_{\bf n}|\widetilde{\bf w}_{\bf n}|\,{\rm e}^{{\rm i}\bf n\cdot x}\right)
\left(\sum_{\bf n}|\widetilde{\bf w}_{\bf n}||{\bf n}|\,{\rm e}^{{\rm i}\bf n\cdot x}\right)
\left(\sum_{\bf n}|\widetilde{\bf w}_{\bf n}||{\bf n}|^{2\gamma}\,{\rm e}^{\rm-i\bf n\cdot x}\right)
\,d{\bf x}$$
$$\le(2\pi)^{-3}{|{\bf w}|}_{6/(3-2\gamma)}
\left|\sum_{\bf n}|\widetilde{\bf w}_{\bf n}||{\bf n}|\,{\rm e}^{{\rm i}\bf n\cdot x}\right|_3
\left|\sum_{\bf n}|\widetilde{\bf w}_{\bf n}||{\bf n}|^{2\gamma}\,{\rm e}^{\rm-i\bf n\cdot x}\right|_{6/(1+2\gamma)}$$
$$\le(2\pi)^{-3}C_\gamma\|{\bf w}\|_\gamma C_{1/2}\|{\bf w}\|_{3/2}C_{1-\gamma}\|{\bf w}\|_{1+\gamma}
\le(2\pi)^{-3}C_\gamma C_{1/2}C_{1-\gamma}\|{\bf w}\|_\gamma^{1/2+\gamma}
\|{\bf w}\|_{1+\gamma}^{5/2-\gamma}$$
\begin{equation}
\le(Q_2/2)\|{\bf w}\|_\gamma^{2(1+2\gamma)/(2\gamma-1)}+\eta_3\|{\bf w}\|_{1+\gamma}^2,
\label{t3}\end{equation}
where we have denoted
$$Q_2=2(2\gamma-1)\left({5-2\gamma\over\eta_3}\right)^{(5-2\gamma)/(2\gamma-1)}
\left({C_\gamma C_{1/2}C_{1-\gamma}\over4(2\pi)^3}\right)^{4/(2\gamma-1)}.$$

Relations \rf{last}--\rf{t3} imply the energy-type inequality
$${d\over dt}(\|{\bf w}\|^2_\gamma+\alpha^2\|{\bf w}\|^2_{\gamma+s})
\le(2\eta_1)^{-1}\beta^2\|{\bf w}\|^2_\gamma+2Q_1\alpha^2\|{\bf w}\|^2_{\gamma+s}
+Q_2\|{\bf w}\|_\gamma^{2(1+2\gamma)/(2\gamma-1)}.$$
Therefore,
$${d\xi\over dt}\le q\xi+Q_2\xi^{(1+2\gamma)/(2\gamma-1)},$$
where
\begin{equation}
\xi(t)\equiv\|{\bf w}\|^2_\gamma+\alpha^2\|{\bf w}\|^2_{1/2+\gamma},\qquad
q\equiv\max((2\eta_1)^{-1}\beta^2,2Q_1).
\label{xiFT}\end{equation}
Integrating this inequality, we obtain a bound
\begin{equation}
\xi(t)\le{\rm e}^{qt}\left((\xi(0))^{-(\gamma-1/2)^{-1}}-(Q_2/q)
({\rm e}^{qt(\gamma-1/2)^{-1}}-1)\right)^{-(\gamma-1/2)}\equiv\varphi(t),
\label{FTNSV}\end{equation}
valid for
\begin{equation}
t<t_*\equiv{2\gamma-1\over2q}\ln\left(1+{q\over Q_2(\xi(0))^{1/(\gamma-1/2)}}\right).
\label{star}\end{equation}
From transformation \rf{TFtr}, relation \rf{xiFT} and this bound we infer

{\it Theorem 5.} Suppose $0<s\le1/2<\gamma\le1$ and the initial
condition of a solution to the Voigt-type regularisation of the force-free
three-dimensional Navier--Stokes equation \rf{VNSeq} belongs to the Sobolev
space $H_{\gamma+1/2}(T^3)$. Then for $t<t_*$ the solution satisfies the bound
$$\Gn{\bf v}({\bf x},t)\Gn_{\beta t,\,\gamma}\le\sqrt{\varphi(t)}.$$
Here a positive constant $\beta$ satisfies inequality \rf{xbeta}, if $s=1/2$,
and is arbitrary otherwise, $\varphi$ and $t_*$ are defined by formulae
\rf{FTNSV} and \rf{star}, respectively, where $\xi(0)$ is determined applying
\rf{xiFT} to ${\bf w}({\bf x},0)={\bf v}^{\rm(in)}$.

This shows that analyticity of solutions to the Voigt-type regularisation
of the Navier--Stokes equation \rf{VNSeq} for $s=1/2$ emerges instantaneously,
provided the initial conditions are in the Sobolev space $H_{\gamma+1/2}(T^3)$.
For small values of the regularising parameter $\alpha$ the bounds $\varphi(t)$
and $t_*$ are uniform in $\alpha$ (note that $q$ is independent of sufficiently
small $\alpha$).

\section{Concluding remarks}

We have explored a new approach to derivation of inequalities for
Gevrey--Sobolev norms of solutions to evolutionary partial differential
equations, in which a suitable nonlinear transformation of a solution
in the Fourier space introduces a feedback between the norm of the transformed
solution and the first index of the norm \rf{GeSo}. Several examples
of application of our technique were discussed.

We have proved that if initially a three-dimensional flow is analytic, then
analyticity in spatial variables is preserved by the Euler equation
on the interval $[0,t_*)$ \rf{tstar} --- this is implied by the bound
\rf{vbndf} for Gevrey--Sobolev norms of a solution.
(Alternatively, a bound for a Gevrey--Sobolev norm of the solution
can be obtained \cite{Titi} from finite-time bounds for solutions to the Voigt
regularisation, which are uniform in the small parameter in the regularising
term ($\alpha$ in \rf{Veq}); such uniform bounds can be derived \cite{Titi}
by application of the method similar to the one used in \cite{LaTi}.)

Pivotal to our technique is transformation \rf{subs}, which results
in emergence of a new term,
$$\beta\varepsilon|{\bf n}|\,\|{\bf w}({\bf x},t)\|_{s+3/2}
^{-1-\varepsilon}\,\widehat{\bf w}_{\bf n}\,{d\over dt}\|{\bf w}({\bf x},t)\|_{s+3/2},$$
in the equations \rf{sEharm} and \rf{sBharm} governing the evolution
of a solution (more precisely, of its transformed Fourier coefficients) to
the Euler and Burgers equations. This term represents a new pseudodifferential
mildly diffusive operator, controlling the norms in $H_{s+3/2}$; it is
analogous to the term
$$-\dot{\tau}(t)\|(-\nabla^2)^{r+1/2}{\bf w}\|^2$$
appearing in literature on analyticity of solutions to PDE's of hydrodynamic
type \cite{Kuka,LeOl}.
(Our assumption $s\le 1/2$ is technical, similar bounds can be obtained
in other Sobolev spaces for larger $s$ in a slightly different way.) The order
of this operator cannot be increased without acquiring difficulties in
deriving bounds for the exponent in the energy balance equation \rf{enin}.
Evidently, by the same construction one can derive an identical bound on the
interval $[0,t_*)$ for solutions to the Navier--Stokes equation. This bound is
uniform in viscosity and can be employed for a study of convergence
of solutions to the Navier--Stokes equation to the solution to the Euler
equation when viscosity tends to zero (such a study is beyond the scope
of the present paper).

Solutions to the Burgers equation are known to develop shocks at finite time
\cite{Four} due to intersection of characteristics; this rules out global
(in time) analyticity in spatial variables of its solutions. That our bounds
for the Gevrey--Sobolev norms of solutions to the Euler and Burgers equations
are identical suggests, that the bounds are rough. Other nonlinear
transformations of solutions can probably yield more accurate bounds. A clear
drawback of our technique is associated with its relative simplicity: we cannot
use it to demonstrate persistence of analyticity of solutions to the Euler
equation, i.e. the fact that the solution is analytic on any interval
$0\le t<t_c$, where it is continuously differentiable \cite{Bard2} (persistence
of non-analytic Gevrey-class regularity of space-periodic solutions to the Euler
equation is shown in \cite{Kuka}).

We have also investigated the Voigt-type regularisations \rf{Veq} of the Euler
equation, involving the regularising term $\alpha^2\partial/\partial t\,(-\nabla^2)^s\bf v$.
and proved boundedness of Gevrey--Sobolev norms of
solutions to the regularised Euler equation for $s>5/6$. Thus, a regularisation,
milder than the one investigated in \cite{LaTi} (for $s=1$), suffices
to guarantee global regularity and analyticity of the solutions. Our technique
enables us to decrease the order of nonlinearity of the energy-type inequality
for the norm of the transformed solution; as a result, our bounds are global
in time and exhibit a polynomial or, at most, exponential growth in time.

Finally, we have explored the Voigt-type regularisation \rf{VNSeq} of the
Navier--Stokes equation with the regularising term for $s>1/2$, which is milder
than the regularisation considered in \cite{KaLeTi} ($s=1$) and the
regularisation of the Euler equation ($s>5/6$) studied here. We have obtained
global exponential in time bounds for Gevrey--Sobolev norms. For $s=1/2$, both
the method \cite{FoTe} and our technique can be applied. Consequently, if
the initial velocity is in the Sobolev
space $H_{\gamma+1/2}(T^3)$, then the solution has a bounded Gevrey--Sobolev
norm \rf{FTNSV} on the time interval $[0,t_*)$ defined by \rf{star}. Furthermore,
starting at any point $t_0$, $0<t_0<t_*$, we obtain a global (in time) bound
\rf{hNSVGS} (performing optimisation of the bounds in $t_*$ and other parameters
is desirable). Unlike all other bounds for Gevrey--Sobolev bounds
that we have derived in this paper, the bound for solutions to the Voigt-type
regularisation \rf{VNSeq} of the Navier--Stokes equation for $s=1/2$ does not
deteriorate in time --- neither the indices of the norm, nor the r.h.s. depend
on time in inequality \rf{hNSVGS}.

\section*{Acknowledgments}

This work was triggered by discussions with Uriel Frisch while walking
around the pine tree park on the premises of the Observatory of Nice. I have
benefited from fruitful discussions with Edriss Titi. My research visits to
the Observatoire de la C\^ote d'Azur in the autumns 2009 and 2010 were supported
by the French Ministry of Education. My work was also partially financed by the
grants ANR-07-BLAN-0235 OTARIE from Agence nationale de la recherche, France,
and 07-01-92217-CNRSL{\_}a from the Russian foundation for basic research.

\end{document}